\input amstex
 \documentstyle{amsppt}
 \UseAMSsymbols
\input cyracc.def

\catcode`\@=11
\font@\tencyr=wncyr10
\font@\eightcyr=wncyr8
\font@\fivecyr=wncyr5
 \font@\sixcyr=wncyr6
 \font@\sevencyr=wncyr7

  \font@\tencyss=wncyss10
\font@\tencyi=wncyi10
 \font@\tencysc=wncysc10

 \font@\eightcyi=wncyi8
\catcode`\@=13

\addto\tenpoint{}
\addto\tenpoint{}

 \addto\eightpoint{}
     \addto\tenpoint{}
 \addto\eightpoint{}

\catcode`\@=11
\font@\tencyb=wncyb10
\font@\eightcyb=wncyb8

\catcode`\@=13
\addto\tenpoint{}

 \addto\tenpoint{}
\NoBlackBoxes

\document

\UseAMSsymbols

   \def\mz{\cprime}

   \def\a{\overline{\bold a}}

    \def\b{\overline{\bold b}}

    \def\c{\overline{\bold c}}

    \def\le{\leqslant}

    \def\ge{\geqslant}

    \def\B{{\Cal B}}

\def\intr{\frac 1 {2\pi}\int_{-\infty}^\infty}

 \def\inti{\int_{-i\infty}^{i\infty}}

\def\pii{\frac 1 {2\pi i}}

\def\NR{Nassrallah--Rahman }

\def\F{{}_2F_1}

\def\FF{{}_3F_2}

\def\FFF{{}_4F_3}

\def\R{{\Bbb R}}

\def\epsilon{\varepsilon}

\def\Re{{\roman Re\,}}

\def\sc{\sl}

\topmatter
\title
 Beta-integrals and finite orthogonal systems
of Wilson polynomials
\endtitle
\author
Neretin Yu.A.
\endauthor 
\thanks Supported by the grant {\rm NWO} 047-008-009.
This work was partially done during my visits to
the Institute of Mathematical Researches (Berkeley)
and the Erwin Schr\"odinger Institute of Mathematical Physics
(Vienna).
 I thank the administrators of these institutes for their hospitality.
       \endthanks
\endtopmatter

 In this paper, we derive the beta-integral
$$ \frac 1 {2\pi}
\int_{-\infty}^\infty
\left|\frac{\prod_{k=1}^3
            \Gamma(a_k+is)}
     {\Gamma(2is)\Gamma(b+is)}\right|^2 ds=
\frac{\Gamma(b-a_1-a_2-a_3)
\prod_{1\le k< l\le 3}\Gamma(a_k+a_l)}
 {\prod_{k=1}^3 \Gamma(b-a_k)}\tag 0.1
$$
(I could not find this integral in literature).
We  construct the system of orthogonal polynomials
related to this integral and also systems
related to the Askey integral
  (0.5)
and the ${}_5\,H_5$--Dougall formula (0.6).
It turns out that all these systems are the Wilson polynomials
(see \cite{Wil}, \cite{AAR}, \cite{KS})
outside the domain of positivity of the classical weight.

\smallskip

{\bf 0.1.} {\bf Beta-integrals.}
The integral (0.1) is a representative of a large family
of so-called beta-integrals.
We mention only the beta-integrals that
in either case will appear
in this paper:
{\it the de Branges--Wilson integral} \cite{dB},
\cite{Wil} (see also \cite{AAR}, 3.6)
$$
\intr
\left|\frac{\prod_{k=1}^4 \Gamma(a_k+is)}
  {\Gamma(2is)}\right|^2 ds=
\frac{\prod_{1\le k< l\le 4}\Gamma(a_k+a_l)}
     {\Gamma(a_1+a_2+a_3+a_4)},\tag 0.2
$$
{\it the Second Barnes Lemma}
(see \cite{AAR}, Theorem 2.4.3)
$$
\multline
\intr
\frac{\Gamma(a_1-is) \Gamma(a_2-is)
       \Gamma(b_1+is) \Gamma(b_2+is) \Gamma(b_3+is)}
  {\Gamma(a_1+a_2+b_1+b_2+b_3+is)}ds=
 \\
=\frac{\prod_{k=1}^3 \Gamma(a_1+b_k)  \Gamma(a_2+b_k)}
 { \prod_{1\le k< l\le 3}\Gamma(a_1+a_2+b_k+b_l)},
\endmultline   \tag 0.3
$$
{\it the Nassrallah--Rahman integral} (see \cite{Ask2},
\cite{RS})
$$
\intr \left| \frac{\prod_{j=1}^5 \Gamma(a_j+is)}
   {\Gamma(2is)\Gamma(\sum_1^5 a_j+is)}\right|^2ds=
2\frac{\prod_{1\le k<l \le 5}\ \Gamma(a_k+a_l)}
  {\prod_{k=1}^5\Gamma(a_1+a_2+a_3+a_4+a_5-a_k)}
,\tag 0.4
$$
and  {\it the Askey integral} (see \cite{RS})
$$
\intr \frac{\Gamma(1-2s)\Gamma(1+2s)}
   {\prod_{k=1}^4\Gamma(a_k+s)\Gamma(a_k-s)}ds
=
\frac{\Gamma(a_1+a_2+a_3+a_4-3)}
     {\prod_{1\le k< l \le 4} \Gamma(a_k+a_l-1)}.
\tag 0.5
$$
In the last case, the integrand has  poles
on the integration contour
 at  the points
$s=\pm n/2, n=1,2,\dots$.
We can understand this integral in two ways.

The first way.
Fix $\alpha$ such that
$0<\alpha<1$.
 Consider the integral
from $-M-\alpha$ to $ N+\alpha$, where $M,N$
are positive integers.
We  understand this integral as the principal
value near each pole.
 For $\Re\sum a_j<3$,
the limit of this integral as
$M,N\to+\infty$ exists.

The second way. Consider the integral
from
  $-M-\alpha$ to $ N+\alpha$
over the contour that passes the poles from above
and consider the limit as  $M,N\to+\infty$.
The both variants coincide since the integrand is an even
function
 (the corrections in the pairs
 $\pm k/2$ of poles cancel).

More complete lists of beta-integrals are contained in
 \cite{Ask2},
\cite{RS},
there are also numerous
 $q$-analogs
(see \cite{AW}, \cite{Ask2}, \cite{RS}),
discrete analogs and multivariate analogs
(see \cite{Gus1}, \cite{Gus2})
 of beta-integrals.

Among discrete analogs,
we need  the
{\it Dougall formula}
$$
\sum_{n=-\infty}^\infty
\frac{\alpha+n}
{\prod_{j=1}^4 \Gamma(a_j+\alpha+n)  \Gamma(a_j-\alpha-n)}
=\frac{\sin 2\pi\alpha}{2\pi}
  \frac{\Gamma(a_1+a_2+a_3+a_4-3)}
   {\prod_{1\le j<k\le 4}\Gamma(a_j+a_k-1)}
.\tag 0.6
$$

{\sc Remark.} In formulas (0.1)--(0.4),
we assume that the parameters
 $a_j$, $b_j$  are real and positive.
To be definite,
consider   integral (0.2).
Represent it
in the form
$$
\intr
\frac{\prod_{k=1}^4 \Gamma(a_k+is)\Gamma(a_k-is)}
  {\Gamma(2is)\Gamma(-2is)} ds=
\frac{\prod_{1\le k< l\le 4}\Gamma(a_k+a_l)}
     {\Gamma(a_1+a_2+a_3+a_4)},
$$
By the analytic continuation, this equality
is valid for  $\Re a_j>0$.
We can  replace this restriction by a weaken variant
 $a_k+a_l\ne 0,-1,-2,\dots$
(for all the pairs $k,l$), but in this case we must change the contour
of  integration
 (see Subsection 0.4)
or  add  the sum of residues in the right-hand side
of the equation.
 We can understand all the integrals (0.1)--(0.4)
and (2.2)--(2.5) in a similar way.

\smallskip

{\bf 0.2.} {\bf  Nassrallah--Rahman integral.}
Our derivation of the integral
 (0.1) is very simple, but the integral
itself can be obtained by degenerations
and hypergeometric transformations
from the following wonderful
 \NR    integral
$$
\int_{-\infty}^\infty
\left|                     \frac{
\Gamma(a+is)\Gamma(b+is)\Gamma(c+is)\Gamma(d+is)
   \Gamma(f+is)\Gamma(g+is)}
 {\Gamma(2is)\Gamma(\lambda+is)\Gamma(\mu+is)}
\right|^2 ds =  \tag 0.7
$$
$$
=\frac{4\pi^3 \sin(a+b+c+d)\pi
  \Gamma(a+b)\Gamma(a+c)\Gamma(a+d)\Gamma(1+2a)}
  {\sin(b+c)\pi\,\sin(b+d)\pi\, \sin(c+d)\pi}
\times
$$
$$
\times
\frac{
 \Gamma(a+f)\Gamma(f-a)\Gamma(g+a)\Gamma(g-a)}
 {\Gamma(1+a-b) \Gamma(1+a-c)  \Gamma(1+a-d)
\Gamma(\lambda+a)\Gamma(\lambda-a)\Gamma(\mu+a)\Gamma(\mu-a)}
\times$$
$$\times
{}_9 F_8\left[
\matrix
2a,  1+a, a+b,  a+c,  a+d,  a+f, a+g, 1+a-\lambda,
     1+a-\mu\\
a,1+a-b,1+a-c,1+a-d,1+a-f,1+a-g,a+\lambda, a+\mu
\endmatrix; 1
\right]
$$
$$
+\dots,
$$
where $\cdots$ denotes the sum of the similar summands,
where
$a$ changes by  $b,c,d$.

We discuss several integrals (2.2)--(2.5)  intermediate between
 (0.2), (0.1)
and
the \NR integral.
All these integrals can be obtained
from  (0.7) by degenerations and hypergeometric transforms,
however I never have seen the final formulas
(2.2)--(2.5).   Nevertheless they deserve to be written
(especially the integral representation (2.2)
for $\FF(1)$), furthermore, our calculations
are very simple conceptually and technically.

\smallskip

{\bf 0.3.} {\bf Finite systems of orthogonal polynomials.}
We write explicitly the system
of polynomials $p_n(s^2)$  orthogonal with respect to the weight
$$w_1(s)=\left|\frac{\prod_{k=1}^3
            \Gamma(a_k+is)}
     {\Gamma(2is)\Gamma(b+is)}\right|^2 ds.$$
Observe that this weight decreases as
$$s^{-2(b-a_1-a_2-a_3)-1},$$
hence only a finite number of the moments
$$\int_0^\infty s^{2k} w_1(s)ds$$
exists. Thus our system of orthogonal polynomials
is finite by the definition.

We also construct the system of polynomials orthogonal with respect to
the discrete weight
$$w_2(s)=\sum c_n\delta_n,$$
where $\delta_n$ is the  delta-function supported by
the point  $n$,
and $c_n$ are the summands of the Dougall formula
 (0.6).

This weight also polynomially decreases.
For some $a_j$ and $\alpha$   the weight is positive, in some cases
it is sign indefinite.
The latter is not very important for play formulas.

Finally, we construct a system of polynomials
(again, this system is finite)
associated with the weight  $w_3(s)$,
where $w_3(s)$ is the integrand in the Askey integral
 (0.5).
This weight is not positive; moreover, it is not
a measure, but we can integrate with respect
to  it.

It turns out that all the three systems
consist of the Wilson
polynomials.
Furthermore, the systems related to the weights
$w_2(s)$ and $w_3(s)$ coincide
(up to a correction factor in the orthogonality
relations).

In first time, finite systems of orthogonal polynomials
(with respect to a continuous weight) were
discovered by
Romanovski \cite{Rom} in 1929
(he constructed analogues of  the Jacobi
polynomials).
Askey  \cite{Ask1} obtained the system of polynomials
related to the Ramanujan
 beta-integral. Various finite families of
orthogonal polynomials
were considered in a series of papers of
Lesky \cite{Les1}--\cite{Les4}.
Some applications of such constructions are contained in
  \cite{Pee}, \cite{BO}.

\medskip

{\bf 0.4.} {\bf Notation.}
$(a)_k:=a(a+1)\dots(a+k-1)$ is the Pochhammer symbol.
$${}_q\,F_{p}
\left[\matrix a_1,\dots a_q\\b_1,\dots,b_p\endmatrix
;z \right]:=\sum_{k=0}^\infty
\frac{ (a_1)_k\dots (a_q)_k }     {k!\,(b_1)_k\dots (b_p)_k} z^k
$$
is the generalized hypergeometric function.
Next,
$$
\Gamma\left[\matrix u_1,\dots u_l\\
       v_1,\dots,v_m\endmatrix\right]:= \frac
  {\Gamma(u_1)\dots\Gamma(u_l)}
   {\Gamma(v_1)\dots\Gamma(v_m)}.
$$

  Mellin--Barnes integrals are defined by
$$   \pii
\inti
 \Gamma\left[\matrix a_1+s,\dots,a_m+s,
 b_1-s,\dots,b_n-s
  \\ c_1-s,\dots,c_k-s, d_1+s,\dots,d_l+s
\endmatrix\right] z^{-s}ds.
$$
We suppose that the integration is given
over the contour
passing from $-i\infty$ to $+i\infty$
and separating the left series of the poles
$$s=-a_1-j,\dots, -a_m-j,\qquad
\text{where}\quad
j=0,1,2,\dots,\infty,$$
from the right series of the poles
$$s=b_1+j,\dots, s=b_n+j.$$
We assume that none of  poles from the left series
coincides with a pole from the right series
(this implies the existence of a contour).
Using the standard rules, we can represent
a Mellin--Barnes  integral
 as a finite sum of hypergeometric functions
with  $\Gamma$-factors.
This rules
(for $m+n\ge k+l$,
bellow this is always satisfied)
are contained in the books
of Slater
\cite{Sla}
and Marichev \cite{Mar}, see also
\cite{PBM}.

\medskip

{\bf \S 1.} {\bf  Preliminaries and preliminary
calculations}

\medskip

{\bf 1.1.} {\bf Index hypergeometric transform.}
Fix $a,b>0$. Let $f$ be a function on the half-line
 $x\ge 0$.
The {\it index hypergeometric transform}
$J_{b,c}$ is defined by the formula
$$
J f(s)=\frac 1{\Gamma(a+b)}
\int_0^\infty f(x) \F(a+is, a-is; a+b; -x)
x^{a+b-1}(1+x)^{a-b}dx
.\tag 1.1
$$
The inverse transform is given by
$$
J^{-1} g(x)=\frac 1{\pi\Gamma(a+b)}
\int_0^\infty g(s) \F(a+is, a-is; a+b; -x)
\left|\frac{\Gamma(a+is)\Gamma(b+is)}{\Gamma(2is)}\right|^2
ds
.\tag1.2
$$
The transformation
 $J$ is a unitary operator
$$
L^2\Bigl(\R_+,x^{a+b-1}(1+x)^{a-b}\Bigr)\to
L^2\Bigl(\R_+,
   \left|\frac{\Gamma(a+is)\Gamma(b+is)}{\Gamma(2is)}\right|^2
\Bigr)
$$
The unitarity condition
({\it the Plancherel formula})
in the explicit form is
$$
\int_0^\infty f_1(x)\overline{ f_2(x)}x^{a+b-1}(1+x)^{a-b}dx=
\frac1\pi\int_0^\infty g_1(s)\overline {g_2(s)}
 \left|\frac{\Gamma(a+is)\Gamma(b+is)}{\Gamma(2is)}\right|^2
ds
\tag 1.3
$$
On properties of this transform, see
\cite{FJK}, \cite{Koo1}--\cite{Koo3},
\cite{Yak}, \cite{Ner}.
The operator $J$ is called also
{\it the Olevsky transform, the  Jacobi transform,
the Fourier--Jacobi transform and the generalized
Fourier transform}, it was discovered by Hermann Weyl
in 1910.

\smallskip

{\bf 1.2.} {\bf  Gauss hypergeometric function.}
We use some simple properties of the hypergeometric functions
 $\F$,
in particular:
the {\it Gauss formula} (see \cite{HTF1}, (2.1.14),
\cite{AAR}, Theorem 2.2,2)
$$
\F(a,b;c;1)=\frac{\Gamma(c)\Gamma(c-a-b)}
{\Gamma(c-a)\Gamma(c-b)};\qquad \Re(c-a-b)>0
.\tag 1.4
$$
and the {\it Bolza  formula} (see \cite{HTF1}
(2.1.22-23))
$$\gather
\F(a,b;c;-x)=(1+x)^{-a}\F\bigl(a,c-b;c;\frac x{x+1}\bigr)=
\tag 1.5\\=
(1+x)^{-b}\F\bigl(c-a,b;c;\frac x{x+1}\bigr)=
(1+x)^{c-a-b}\F\bigl(c-a,c-b;c;-x\bigr)
.\tag 1.6
\endgather
$$
Recall that the function
 $\F(a,b;c;x)$ admits the analytical continuation
to the half-line $x<0$.
The asymptotics of the function
$\F(a+is,a-is;a+b;-x)$ as $x\to+\infty $
has the form
$$\F(a+is,a-is; a+b; -x)=\psi(x) x^{-a},$$
where $\psi$  is a bounded function
for $s\ne 0$ and $\psi(x)/\ln(x)$ is a bounded function
for $s=0$
(for a more explicit expression, see
 \cite{HTF1},(2.3.2.9)).

The asymptotics of the function
$\F(a+is,a-is; a+b; -x)$
as $s\to +\infty$ has the form
$$\F(a+is,a-is; a+b; -x)=\lambda(s) s^{-a-b+1/2}.$$
where $\lambda(s)$ is a bounded function
(see the Watson formula  \cite{HTF1}, (2.3.2.17)).

Recall three representations of the Gauss hypergeometric function
as Barnes integrals
see \cite{HTF1}, (2.1.15), \cite{Mar},
\S 2.10, formulas 11(1), 19(1), 22(1),
\cite{AAR}, Theorems 2.4.1, 2.4.2, \cite{PBM},
8.4.49(13), (20), (22);
below
  $x>0$

$$ \pii
\frac{\Gamma(c)}{\Gamma(a)\Gamma(b)}
\inti
\frac{\Gamma(a-s)\Gamma(b-s)\Gamma(s)}
{\Gamma(c-s)} x^{-s}ds =
\F(a,b;c;-x);          \tag 1.7
$$
$$\multline
\pii \frac{\Gamma(c)}
{\Gamma(a)\Gamma(b)\Gamma(c-a)\Gamma(c-b)}
\inti
\Gamma(s)\Gamma(s+c-a-b)\Gamma(a-s)\Gamma(b-s)x^{-s}\,ds =\\
=
\F(a,b;c;1-x)         ;
\endmultline\tag 1.8
$$
$$    \multline
\pii
\Gamma(c)\inti
\frac{\Gamma(s)\Gamma(s+c-a-b)}
     {\Gamma(s+c-a)\Gamma(s+c-b)}x^{-s}ds
=\\=
\left\{ \matrix(1-x)^{c-1}
       \F(a,b;c;1-x); \quad x> 1,\\
     0; \quad x< 1\endmatrix\right.
\endmultline \tag 1.9
$$

{\bf 1.3.} {\bf $\Gamma$-function.}
We shall use the following integral representation
of the beta-function
$$
\int_0^\infty \frac{x^{\alpha-1}dx}
{(1+x)^\rho}=B(\alpha,\rho-\alpha),
$$
The substitution
 $x/(x+1)=z$ transforms this integral to
the usual definition of the beta-function.

Asymptotics of the
 function $\Gamma(z)$  is given by the
 Stirling formula
$$\Gamma(z)=\sqrt{2\pi}
z^{z-1/2}e^{-z}(1+O(1/z)); \qquad z\to\infty
\qquad \arg z<\pi-\epsilon.
$$
In particular (see \cite{HTF1}, (1.18.3), (1.18.6)),
$$\frac{\Gamma(z+a)}{\Gamma(z+b)}\sim
z^{a-b},\qquad z\to\infty,  \qquad \arg z<\pi-\epsilon,
$$
$$
|\Gamma(a+is)|\sim\sqrt{2\pi} |s|^{a-1/2}e^{-|s|\pi/2},
\qquad s\to +\infty.
$$

{\bf 1.4. \bf  Mellin transform.}
Consider a function  $f$ on the half-line   $x\ge 0$.
Its Mellin transform is defined by the formula
$$
f^*(s)=\int_0^\infty x^s f(x)\frac{dx} x .
$$
The inversion formula is
$$
f(x)=\pii\inti
f^*(s) x^{-s}ds                  .
$$

The Mellin transform takes the convolution
$$
\int_0^\infty f_1(x/t)f_2(t)\frac {dt} t
$$
to the product of the Mellin transforms
$f_1^*(s)f_2^*(s)$.
For precise statements, see for instance
\cite{Mar}.

Let $a>0$. Obviously, the function $f(ax)$
is mapped to
 $a^{-s} f^*(s)$.
Also, the function $x^\alpha f(x)$
corresponds to $f^*(s+\alpha)$.

{\bf 1.5.} We need  some integrals of the form
$$               \multline
\int_0^\infty
\frac{x^{\alpha-1}}{(x+z)^\rho} \F(p,q;r;-x)\,dx=
\\=
\pii  z^{\alpha-\rho}
\Gamma\left[\matrix r\\ p,q,\rho
\endmatrix\right]
\inti\Gamma\left[
\matrix s+\alpha,\rho-s-\alpha,p+s,q+s,-s\\
r+s\endmatrix\right]
z^s    ds
\endmultline\tag 1.10  $$
 The Mellin--Barnes
integral in the right-hand side
can be expressed as a linear combination of two
 $\FF$.

As a quite standard (see \cite{Mar})  pattern,
we prove the identity
 (1.10); this statement also is quite standard
(see
 \cite{PBM}, formula (2.21.1.16)).
For this, represent
$x^{\alpha}(x+z)^{-\rho}$ as the inverse Mellin transform
 $$x^{\alpha}(x+z)^{-\rho}=
\frac{z^{\alpha-\rho}}{2\pi\Gamma(\rho)}\inti
\Gamma(s+\alpha)\Gamma(\rho-s-\alpha)z^s x^{-s} ds
$$
(for a verification, it is sufficient
to evaluate the direct Mellin transform).
Applying the Barnes integral (1.7)
and the convolution formula
for the Mellin transform,
we obtain (1.10).

\smallskip

{\bf 1.6.} {\bf Corollaries from formula (1.10).}

\smallskip

{\sl Lemma 1.1.} {
{\rm  {a)}} \it
The
operator $J$ takes the function
$$
(1+x)^{-a-c}\qquad\ \text{to}\qquad \qquad
\frac{\Gamma(c+is)\Gamma(c-is)}
     {\Gamma(c+a)\Gamma(c+b)}
.$$
{\rm  {b)}}
The operator $J$ takes the function
$$
 \frac{(1+x)^{b-a}}
       {(x+z)^{c+b}}
\qquad\ \text{to}\qquad \qquad
      \frac{\Gamma(c+is)\Gamma(c-is)}
     {\Gamma(c+a)\Gamma(c+b)}
\F\left[\matrix c+is,c-is\\c+a \endmatrix ; 1-z\right]
.$$
{\rm  {c)}}
The  operator $J$ takes the function
$$
 x^{-u-a}   \qquad\ \text{to}\qquad \qquad
\frac{\Gamma(-u+b)}{\Gamma(a+u)}
\cdot
\frac{\Gamma(u+is) \Gamma(u-is)}
      {\Gamma(b+is) \Gamma(b-is)}
.$$
 }

{\sl Proof.}
We must substitute these
 $f$
to formula (1.1).
We can simply refer to the rows
(2.21.1.4), (2.21.1.5),  (2.21.1.15)
of
Prudnikov, Brychkov, Marichev, vol.3,
or to the rows (7.512) of
Ryzhik--Gradshteyn.

But it is better to explain what happened
and why for some values of the parameters the integral
 (1.10) admits an explicit evaluation
and for some values this is impossible.

  1. For $\alpha=r$ in the right-hand side of  (1.10),
two  $\Gamma$-factors in the integrand cancel.
  The remains is a Barnes integral  of the type  (1.8);
this proves  the statement
{\rm b)}.

  2.  Substituting   $z=1$  to the statement
 {\rm a)}, we obtain the claim  {\rm b)}.

  3. Substitute $z=1$ to (1.10). For  $r=p+q+\rho$,
the right-hand side can be evaluated by the Second Barnes Lemma
 (0.3).  This gives  {\rm c)}.
\hfill$\boxtimes$

\smallskip

 {\bf 1.7.} {\bf Integrals with products of hypergeometric
functions.}
A direct application of the convolution formula for the Mellin
transform gives
$$\multline
\int_0^\infty
x^{\alpha-1}
\F\left[\matrix p,q\\r \endmatrix; -\omega x\right]
\F\left[\matrix u,v\\w \endmatrix; -\tilde\omega x\right]
dx=\\=
\pii \omega^{-\alpha}\Gamma \bmatrix r,w\\ u,v,p,q\endbmatrix
\inti \Gamma\bmatrix
\alpha+s,u+s,v+s,p-\alpha-s,q-\alpha-s,-s\\
 r-\alpha-s,w+s\endbmatrix
   \left(\frac{\omega}{\tilde\omega}\right)^{-s}ds
.\endmultline  \tag 1.11
$$
Substitute
$\alpha=w=r$ to this identity.
Then four $\Gamma$-factors
of the integrand in the right-hand side cancel.
Applying the Barnes formula  (1.8),
we obtain
$$ \multline
\int_0^\infty x^{r-1}
\F\left[\matrix p,q\\r \endmatrix; -\omega x\right]
\F\left[\matrix u,v\\r \endmatrix; -\tilde\omega x\right]
dx=\\ =
\omega^{-r+u}\tilde\omega^{-u}
\Gamma\bmatrix r,r,p-r+u,q-r+u,p-r+v,q-r+v\\
 u,v,p,q,p+q+u+v-2r\endbmatrix\times \\ \times
\F \left[\matrix p-r+u,q-r+u\\p+q+u+v-2r\endmatrix
   ; 1- \frac{\omega}{\tilde\omega} \right]
.\endmultline                      \tag 1.12
$$
It is formula (2.21.9.7) from  Tables \cite{PBM}
(where a prime in the second index of
hypergeometric function in the right-hand side is lost).

We also need  the integral
$$\multline
\int_0^1 z^{\mu-1} (1-z)^{\nu-1}
\F \left[\matrix \alpha,\beta\\ \nu   \endmatrix; 1-z   \right]
\F \left[\matrix \phi,\psi\\ \xi   \endmatrix; 1-z   \right]
dz=\\=
\pii
\Gamma\bmatrix \nu,\xi\\ \phi,\psi,\xi-\phi,\xi-\psi\endbmatrix
\times\\ \times
\inti \Gamma\bmatrix \mu+s, \mu+\nu-\alpha-\beta+s,
\phi+s,\psi+s,-s,\xi-\phi-\psi+s\\
\nu+\mu-\alpha+s, \mu+\nu-\beta+s  \endbmatrix ds
.\endmultline \tag 1.13
$$
To deduce this formula, we represent
 $(1-z)^{\nu-1}  \F( \alpha,\beta; \nu;1-z)$
as the integral (1.9), and represent  $\F( \phi,\psi; \xi; 1-z)$
as the integral (1.8).   It remains to apply the convolution formula.

\smallskip

{\bf 1.9.} {\bf Corollaries from the formula  (1.11).}

\smallskip

{\sl Lemma 1.2.}
{\rm a)} {\it The transform $J$
takes
$$ \F\left[\matrix p+b,q+b\\ a+b\endmatrix
   ;- \dfrac xy \right](1+x)^{b-a}$$
 to
$$
\frac{y^{b-q}\Gamma(a+b)}{\Gamma(p+q)\Gamma(p+b)\Gamma(q+b)}
\cdot \frac{\Gamma(p+is)\Gamma(p-is)\Gamma(q+is)\Gamma(q-is)}
{\Gamma(a+is)\Gamma(a-is)}
 \F \left[\matrix p+is,p-is\\p+q  \endmatrix; 1-y \right]
.$$

 {\rm b)}  The transform $J$
takes $ \F\left[\matrix p+b,q+b\\ a+b\endmatrix
   ;-x \right](1+x)^{b-a}$ to
$$
  \frac{\Gamma(a+b)}{\Gamma(p+q)\Gamma(p+b)\Gamma(q+b)}
\cdot \frac{\Gamma(p+is)\Gamma(p-is)\Gamma(q+is)\Gamma(q-is)}
{\Gamma(a+is)\Gamma(a-is)}
.$$
 { \rm c)} The transform  $J$ takes
 $\F\left[\matrix a+c,a+d\\ a+b+c+d\endmatrix;-x\right]$
to      }
$$\frac
{\Gamma(a+b+c+d)\cdot\Gamma(c+is)\Gamma(c-is)
\Gamma(d+is)\Gamma(d-is)}
{\Gamma(a+c)\Gamma(a+d)\Gamma(b+c)\Gamma(b+d)\Gamma(c+d)}
.$$

{\sl Proof.}
The claim {\rm a}) is an immediate corollary
of formula (1.11). The  substitution  $y=1$ to {\rm a)}
gives the
statement  {\rm b)} (of course, it is more pleasant
to substitute  $\omega=\tilde\omega$
directly to  (1.11).

Let us deduce {\rm c)}. Substituting $f=\F[\dots]$
to (1.1) and applying the Bolza  transformation,
we obtain
$$\frac 1{\Gamma(a+b)}\int_0^\infty
x^{a+b-1} \F\left[\matrix b+d, b+c\\a+b+c+d\endmatrix
;-x\right]
 \F\left[\matrix a+is, a-is\\a+b  \endmatrix ;-x\right]
dx
$$
Furthermore, in the integrand in the right-hand side
of (1.12),
two $\Gamma$-factors cancel and we obtain
$$\multline
\pii\Gamma\bmatrix a+b+c+d\\a+is,a-is,b+c,b+d\endbmatrix
\times\\ \times
\inti \Gamma\bmatrix a+is+t, a-is+t, d-a-t, c-a-t,-t\\
             c+d-t\endbmatrix \,dt
\endmultline
$$
It remains to apply the Second Barnes Lemma (0.3).

 \medskip

{\bf \S 2.}  {\bf Degenerate cases of
the \NR integral}

\medskip

In this section, we apply the Plancherel formula
 (1.3) for the index hypergeometric transform
to various pairs of functions
from Lemmas
 1.1--1.2.

To be brief, we assume
that the parameters $a,b,c,d,e,f, u, v$
are real positive.

{\bf 2.1.} {\bf De Branges--Wilson integral.}
Apply the Plancherel formula
  (1.3) to the pair of functions
$(1+x)^{-a-c}$ and $(1+x)^{-a-d}$,
see Lemma 1.1. In the left-hand side we
have the beta-integral

$$
\int_0^\infty \frac{x^{a+b-1}dx}{(1+x)^{a+b+c+d}} \tag 2.1
.$$
In the right-hand side, we obtain the de Branges--Wilson integral
(up to $\Gamma$-factors).
Evaluating  (2.1), we obtain  (0.2).

This proof of the de Branges--Wilson integral is known,
see the work
 of Koornwinder \cite{Koo2}.

\smallskip

{\bf 2.2.} {\bf The integral  (0.1).}
It is sufficient to apply the Plancherel formula
to the pair of functions
$x^{-u-a}$ and $x^{-v-a}$, see Lemma   1.1.{\rm c}.

\smallskip

{\bf 2.3.}
Apply the Plancherel formula  (1.3)
to the pair of functions
$$
(1+x)^{-a-e}\qquad\text{and} \qquad
\F\left[\matrix a+c, a+d\\a+b+c+d\endmatrix;-x\right]
.$$
In the left-hand side we have
$$\int_0^\infty x^{a+b-1}(1+x)^{-b-e}
 \F\left[\matrix a+c, a+d\\a+b+c+d\endmatrix;-x\right] dx
.$$
We transform this by the Bolza  formula and obtain
$$\int_0^1 z^{a+b-1}(1-z)^{c+e-1}
 \F\left[\matrix a+c, b+c\\a+b+c+d\endmatrix;z\right]
dz
.$$
This is a standard integral representation
for
 $\FF$ (see \cite{AAR}, (2.2.4)),
and finally the left-hand side is
$$
\frac{\Gamma(a+b)\Gamma(c+e)}{\Gamma(a+b+c+e)}
\FF\left[\matrix a+c, b+c, a+b\\ a+b+c+d, a+b+c+e\endmatrix
;1\right]
$$
Equating the left-hand and right-hand sides of the Plancherel formula,
we obtain
$$\multline
\frac{1}{\pi}\int_0^\infty
\left|\frac{\Gamma(a+is)\Gamma(b+is)\Gamma(c+is)\Gamma(d+is)
\Gamma(e+is)}
{\Gamma(2is)}\right|^2ds=\\=
\frac{\Gamma(a+b)\Gamma(a+c)\Gamma(a+d)\Gamma(a+e)
\Gamma(b+c)\Gamma(b+d)\Gamma(b+e)\Gamma(c+d)\Gamma(c+e)}
{\Gamma(a+b+c+d)\Gamma(a+b+c+e)}
\times\\ \times
 \FF\left[\matrix a+c, b+c, a+b\\ a+b+c+d, a+b+c+e\endmatrix
;1\right]
\endmultline \tag 2.2
.$$
Among the $\Gamma$-factors in the numerator,
the factor
$\Gamma(e+d)$ is absent.

{\sl Remark.}
Obviously, the left-hand side of the identity
is symmetric in
 $a,b,c,d,e$. Symmetry of the right-hand side
is equivalent to the
 Kummer formula, see
\cite{AAR}, Corollary 3.3.5.

\smallskip

{\bf 2.4.}
Now apply the Plancherel formula
(1.3) to the pair of functions
$$
  \F\left[\matrix a+c, a+d\\a+b+c+d\endmatrix;-x\right]
\qquad\text{ and} \qquad
   \F\left[\matrix a+e, a+f\\a+b+e+f\endmatrix;-x\right],
$$
see Lemma 1.2.{\rm c}.
In the left-hand side, we have
$$ \multline
\int_0^\infty
x^{a+b-1}
   \F\left[\matrix a+c, a+d\\a+b+c+d\endmatrix;-x\right]
   \F\left[\matrix a+e, a+f\\a+b+e+f\endmatrix;-x\right]
(1+x)^{a-b}dx
=\\=
\int_0^\infty
x^{a+b-1}
   \F\left[\matrix b+d, b+c\\a+b+c+d\endmatrix;-x\right]
   \F\left[\matrix a+e, a+f\\a+b+e+f\endmatrix;-x\right]dx
.\endmultline
$$
Applying formula (1.11),
we obtain
$$\multline
\frac 1\pi
\int_0^\infty
  \left|\frac{\Gamma(a+is)\Gamma(b+is)\Gamma(c+is)\Gamma(d+is)
\Gamma(e+is)\Gamma(f+is)}
{\Gamma(2is)}\right|^2ds=\\ =
\frac 1{2\pi i}
\Gamma(a+c)\Gamma(a+d)\Gamma(c+d)
\Gamma(b+e)\Gamma(b+f)\Gamma(e+f)\times
\\ \times
\inti
\Gamma\bmatrix a+b+s, a+e+s, a+f+s, d-a-s,c-a-s,-s
      \\
       c+d-s,a+b+e+f+s
       \endbmatrix
ds
.\endmultline\tag 2.3
$$
The expression in the right-hand side can be represented
as a linear combination of three functions
 $\FFF(1)$.

\smallskip

{\bf 2.5.}  Finally, we apply the Plancherel formula
to the pair of functions
$$
\F\left[\matrix p+b,q+b\\a+b\endmatrix;-x\right]
    (1+x)^{b-a}
\quad\text{and} \quad
 \F\left[\matrix u+b,v+b\\a+b\endmatrix;-x\right]
    (1+x)^{b-a}
$$
(their index transforms are evaluated
in Lemma  1.2{\rm b}).

We must evaluate
the integral
$$
\int_0^\infty
x^{a+b-1}
(1+x)^{b-a} \F\left[\matrix p+b,q+b\\a+b\endmatrix;-x\right]
 \F\left[\matrix u+b,v+b\\a+b\endmatrix;-x\right]
dx              .
$$

Using the Bolza  formula, we obtain
$$  \multline
\int_0^1
y^{a+b-1}
(1-y)^{u+p}
 \F\left[\matrix p+b,a-q\\a+b\endmatrix;y\right]
 \F\left[\matrix u+b,a-v\\a+b\endmatrix;y\right]
dy=          \\ =
\int_0^1
(1-z)^{a+b-1}
z^{u+p}
 \F\left[\matrix p+b,a-q\\a+b\endmatrix;1-z\right]
 \F\left[\matrix u+b,a-v\\a+b\endmatrix;1-z\right]
dz               .\endmultline
$$
It remains to apply
 (1.13).
As a result, we get
 $$\multline
\frac 1\pi
\int_0^\infty
\left|\frac{\Gamma(b+is)\Gamma(p+is)\Gamma(q+is)\Gamma(u+is)
\Gamma(v+is)}
{\Gamma(2is)\Gamma(a+is)}\right|^2ds =\\=
\frac 1{2\pi i}
\Gamma\bmatrix u+v,p+q,p+b,q+b\\ a-v,u-v\endbmatrix
\times\\ \times  \inti
 \Gamma\bmatrix u+p+s,u+q+s,b+u+s,a-v+s,v-u-s,-s\\
u+a+s,u+b+p+q+s\endbmatrix      .
\endmultline          \tag 2.4
.$$

{\sl Remark.} The right-hand side can be written in the form
$$
\multline
\Gamma\bmatrix
u+v,p+q,p+b,q+b\\ a-v,u-v
\endbmatrix\times\\ \times
\Biggl\{
\Gamma\bmatrix v-u,u+p,u+q,u+b,a-v\\
  u+a,u+b+p+q\endbmatrix
\FFF\left[\matrix u+p,u+q,u+b,a-v\\
1+u-v,u+a,u+b+p+q\endmatrix;1\right]
+\\+
\Gamma\bmatrix u-v,p+v,q+v,b+v,a-u\\
       v+a,v+b+p+q\endbmatrix
 \FFF\left[\matrix
     p+v,q+v,b+v,a-u\\
     1-v+u,v+a,v+b+p+q
      \endmatrix;1\right]
      \Biggr\}
.\endmultline
$$
 Rahman \cite{Rah}
gives the for   integral (2.4)
a
${}_7 F_6$-expression (see also \cite{GR}, (6.3.11)).
However, it can be easily reduced to the form
 $\FFF$  using the
nonterminating ${}_7F_6$-Whipple transform,
discussed in Bailey's book \cite{Bai}, 4.4, 6.3,  7.5.

\smallskip

{\sl Remark.} The left-hand side of the identity
(2.4)
is symmetric in $b,p,q,u,v$.
In the right-hand side, the symmetry in
 $u,v$ and in $b,p,q$ is obvious.
The transposition of  $b$ and $u$
gives  a 4-term $\FFF$-identity.
Substitute
$$a=v-m,$$
where $m$ is a nonnegative integer, to this identity.
Then two summands vanish
(due $\Gamma(-m)$ in the denominators),
and we obtain  $\FFF$-Whipple transform
(see \cite{AAR}, Theorem 3.3.3;
it is mentioned below  in Subsection
 3.1).

\smallskip

{\bf 2.6.} {\bf  \NR integral (0.4).}
If
$$
a=b+u+v+p+q
,$$
then two $\Gamma$-factors in the right-hand side
of (2.4) cancel.
Applying the Second Barnes Lemma (0.3),
we obtain  (0.4).

\smallskip

{\bf 2.7.} Multiplying   the both sides
of (2.4)
by $\Gamma^{-2}(b)$ and passing to the limit as
 $b\to+\infty$,   we obtain
  $$\multline
\frac
1\pi\int_0^\infty
\left|\frac{\Gamma(p+is)\Gamma(q+is)\Gamma(u+is)
\Gamma(v+is)}
{\Gamma(2is)\Gamma(a+is)}\right|^2ds =\\=
\frac 1{2\pi i}
\Gamma\bmatrix u+v,p+q\\ a-v,u-v\endbmatrix
 \inti
 \Gamma\bmatrix u+p+s,u+q+s,a-v+s,v-u-s,-s\\
u+a+s\endbmatrix ds
\endmultline    \tag 2.5
.$$
(in the right-hand side we have a linear combination
of two functions
 $\FF(1)$).

\medskip

{\bf \S 3. \bf  Finite systems of orthogonal polynomials.}

\medskip

{\bf 3.1.} {\bf Wilson polynomials.}
In the famous work \cite{Wil},
Wilson constructed the polynomials  $p_n(s^2)$
orthogonal with respect to the weight
 $$w(s) =
 \frac 1\pi\left|\frac{\Gamma(a+is)\Gamma(b+is)
 \Gamma(c+is)\Gamma(d+is)}
{\Gamma(2is)}\right|^2;$$
they are defined by the formula
$$\multline
p_n(a,b,c,d;s^2)=\\=(a+b)_n (a+c)_n (a+d)_n
\,\,
\FFF\left[\matrix
-n,n+a+b+c+d-1, a+is,a-is\\
a+b,a+c,a+d\endmatrix; 1\right]
\endmultline    \tag 3.1
.$$
The orthogonality relations have the form
$$        \multline
\int_0^\infty p_n(a,b,c,d;s^2) p_m(a,b,c,d;s^2)
 w(s)ds=\\ =
\frac{n!\Gamma(a+b+n)\Gamma(a+c+n)\Gamma(a+d+n)
  \Gamma(b+c+n)\Gamma(b+d+n)\Gamma(c+d+n)}
  {\Gamma(a+b+c+d+n) (a+b+c+d+2n-1)}       \delta_{m,n}
\endmultline            \tag 3.2
$$
Evidently, the polynomials
$p_n(a,b,c,d;s^2)$ are symmetric with respect to
 $b,c,d$.
The symmetry in four indices
 $a,b,c,d$ is equivalent to
$\FFF$-Whipple transform
mentioned above in
 Subsection 2.5.

\smallskip

{\bf 3.2.} {\bf Proof of orthogonality.}
Several proofs of orthogonality of
the Wilson polynomials are known.
Our purpose is to give a model proof
that works in 3 cases of "exotic orthogonality"
discussed below.

\smallskip

{\sl Lemma 3.1.}
{\it  Let $V$, $W$  be $N$-dimensional linear spaces
with bases  $e_j$, $f_j$ respectively.
Define the sesquilinear  form $\B(\cdot,\cdot)$
on
$V\times W$ by
$$\B(e_k,f_l)=\frac{\Gamma(\nu+k+l)}
             {\Gamma(\mu+\nu+k+l)}
,\tag 3.3$$
where $\mu,\nu$ are fixed. Then the system
of vectors
$$
R_n=\sum_{j=1}^n\frac{(-n)_j(n+\mu+\nu-1)_j}
             {(\nu)_j \,j!} e_j\qquad\qquad
 T_n=\sum_{j=1}^n\frac{(-n)_j(n+\mu+\nu-1)_j}
             {(\nu)_j \,j!} f_j
$$
is biorthogonal, i.e.,
  $\B(R_k,T_l)=0$ for
$k\ne l$. Moreover,  }
$$\B(T_n,R_n)=\frac{n!\,\Gamma(\nu)\Gamma(\mu+n)}
               {\Gamma(\nu+n))\Gamma(n+\mu+\nu-1)
                  (\mu+\nu+2n-1)}
.\tag 3.4
$$

{\sl Proof.} This lemma imitates the orthogonality relations
for the Jacobi polynomials.
Let
 $V=W$ be the space of polynomials
on the segment $[0,1]$ with the scalar product
$$
<p,q>=
\frac 1{\Gamma(\mu)}
\int_0^1 p(x)\overline{q(x)}x^{\nu-1}(1-x)^{\mu-1}dx
.$$
Let $e_n=f_n$ be the function
$x^n$. Then $<e_k,f_l>$ coincides with
(3.3).  Next,
$R_n=T_n$ are the Jacobi polynomials
$P^{\mu-1,\nu-1}(2x-1)$ in the standard notation
(see \cite{HTF2}, 10.8).
 \hfill $\square$

\smallskip

Now let
  $V=W$ be the space
 $L^2$ with respect to the Wilson weight  $w(s)$.
Assume
$$
e_k(s)=\frac{(a+is)_k(a-is)_k}
            {\Gamma(a+c+k)\Gamma(a+d+k)}
,\qquad
f_k(s)= \frac{(b+is)_k(b-is)_k}
            {\Gamma(b+c+k)\Gamma(b+d+k)}
.\tag 3.5
$$
Then by
 (0.2),
$$\multline
\int_0^\infty e_k(s) f_m(s) w(s)\,ds= \\=
\frac 1{\pi\Gamma(a+c+k)\Gamma(a+d+k)
\Gamma(b+c+m)\Gamma(b+d+m)}
\times \\ \times \int_0^\infty
 \left|\frac{\Gamma(a+k+is)
\Gamma(b+m+is)\Gamma(c+is)\Gamma(d+is)}
{\Gamma(2is)}\right|^2ds=\\=
\frac{\Gamma(c+d)\Gamma(a+b+k+m)}
 {\Gamma(a+b+c+d+k+m)}
.\endmultline
$$
We obtain the relations of the type
(3.3) for scalar products.
In notation of  Lemma 3.1, we have
$$
R_n=\alpha_n p_n(a,b,c,d;s^2);\qquad
T_n=\beta_n  p_n(b,a,c,d;s^2)
,$$
where
$\alpha_n, \beta_n$ are normalizing factors
(we omit them to be brief).
It remains to recall that the Wilson polynomials
are symmetric with respect to
 $a,b,c,d$.
Hence $R_n$ and $T_n$ coincide up to a factor
(precisely this place of the proof is surprising).

Formula (3.2) follows from  (3.4).

\smallskip

 {\bf 3.3.} {\bf  Orthogonality relations
associated with the integral (0.1).}
Consider the weight
  $$w_1(s)=
  \frac 1\pi \left|\frac{\Gamma(p+is)\Gamma(u+is)\Gamma(v+is)}
     {\Gamma(2is)\Gamma(q+is)}\right|^2
\tag 3.6
   $$
on the half-line
 $s\ge 0$.

Assume
$$
\gather
e_n(s)=\frac{\Gamma(-u-k+q)}{\Gamma(p+u+k)}
   (u+is)_k(u-is)_k;\\
f_n(s)=\frac{\Gamma(-v-k+q)}{\Gamma(p+v+k)}
   (v+is)_k(v-is)_k
.\endgather
$$
Then
$$ \multline
\int_0^\infty
e_k(s) f_l(s)w_1(s)\,ds=
\frac{\Gamma(u+v+k+l)\Gamma(q-u-v-p-k-l)}
     {\Gamma(q-p)} =  \\
\frac{(-1)^{k+l}\pi}
{\Gamma(q-p)\sin \pi(q-u-v-p)}
\cdot\frac{\Gamma(u+v+k+l)}{\Gamma(1-q+u+v+p+k+l)}
.\endmultline
$$
Observe that the scalar products have the form
 (3.3). Thus we can repeat literally
all the considerations of Subsection 3.2.

There exists a way that is even more simple.

\smallskip

{\sl Lemma 3.2.} {\it  Let
$\ell$  be a linear functional on the space
of even polynomials
$$c_0+c_1 s^2+\dots+c_N s^{2N}$$
and for
$$h_k=(a+is)_k(a-is)_k$$
we have
$$
\ell(h_k)=C\cdot \frac{\Gamma(a+b+k)\Gamma(a+c+k)\Gamma(a+d+k)}
{\Gamma(a+b+c+d+k)}
.\tag 3.7
$$
Then for $k+l\le N$
the Wilson polynomials
$p_k,p_l$
satisfy
$$\ell (p_k(a,b,c,d;\cdot) p_l(a,b,c,d;\cdot))=
\frac {C}{\Gamma(b+c)\Gamma(b+d)\Gamma(c+d)}\cdot
\sigma_k(a,b,c,d)\delta_{k,l},
$$
where $\sigma_n(a,b,c,d)$ is the expression
situated  in the right-hand side of
} (3.2).

\smallskip

In other words, the Wilson polynomials form
an orthogonal system with respect to the scalar product
$$\langle f,g\rangle =\ell(f\overline g)$$
in the space of polynomials.

\smallskip

{\sl Proof.}
Assume
$$\ell(g)=\int_0^\infty g(s)w(s)ds,$$
where $w(s)$ is the Wilson weight.
Then (3.7) is satisfied
and our statement is a rephrasing
of the orthogonality relations for the Wilson polynomials.
\hfill $\square$

\smallskip

Now suppose
$$\ell(g)= \int_0^\infty g(s)w_1(s)ds$$
and
$$h_k(s)=(p+is)_k(p-is)_k.$$
We have
$$
\multline
\ell(h_k)=
\frac{\Gamma(p+v+k)\Gamma(p+u+k)\Gamma(u+v)\Gamma(q-u-v-p-k)}
    {\Gamma(q-p-k)\Gamma(q-u)\Gamma(p-v)} = \\
=\frac{\Gamma(u+v)\sin\pi(p-q)}
{\Gamma(u-q)\Gamma(q-v)\sin\pi(u+v+p-q)}
\cdot \frac{\Gamma(p+v+k)\Gamma(p+u+k)\Gamma(1-q+p+k)}
            {\Gamma (1-q+u+v+p+k)}.
.\endmultline
$$

Therefore, the system of the Wilson polynomials
$$
p_n(p,u,v,1-q;s^2); \qquad 4n<q-p-u-v-1
$$
is orthogonal with respect to the weight (3.6).

\smallskip

 {\bf 3.4. \bf   Orthogonality relations
associated with Dougall formula.}
Consider the weight
$$w_2(t)=\sum_{-\infty}^\infty
\frac{\alpha+t}
{\prod_{j=1}^4 \Gamma(a_j+\alpha+t)  \Gamma(a_j-\alpha-t)}
\delta (t-n)
,$$
where $\delta(s-n)$ is the delta-function
supported by the point $n$.
For uniformity,  substitute
 $t=is-\alpha$ and transform
this expression to the form
$$w_2(s):=\prod_{j=1}^4
 \frac{\sin\pi(a_j+\alpha)}
 {\pi}
\cdot (is)
\sum_{-\infty}^\infty
\biggl\{\prod_{j=1}^4
                 \Gamma(1-a_j-is)\Gamma(1-a_j+is)
   \cdot
\delta (is-\alpha-n)
\biggr\}
.\tag 3.8
$$

To apply
Lemma 3.2, suppose
$$ \ell(g)= \int g(s)w_2(s)ds$$
$$h_k(s)=(1-a_1+is)_k(1-a_1-is)_k.$$
Then
$$\multline
\ell(h_k)=\frac{\sin (2\pi\alpha)}
{2\pi\Gamma(a_2+a_3)\Gamma(a_2+a_4)\Gamma(a_3+a_4)}
\times\\ \times
\frac{\Gamma(a_1+a_2 +a_3+a_4-k-3)}
{\Gamma(a_1+a_2-k-1)\Gamma(a_1+a_3-k-1)\Gamma(a_1+a_4-k-1)}
=\\ =
 \frac{\sin (2\pi\alpha)
\sin\pi(a_1+a_2)\sin\pi(a_1+a_3)\sin\pi(a_1+a_4)}
{2\pi\Gamma(a_2+a_3)\Gamma(a_2+a_4)\Gamma(a_3+a_4)
\sin\pi(a_1+a_2+a_3+a_4)}
\times \\ \times
\frac
{\Gamma(-a_1-a_2+k+2) \Gamma(-a_1-a_3+k+2) \Gamma(-a_1-a_4+k+2)}
{\Gamma(4-a_1-a_2-a_3-a_4+k)}
.\endmultline\tag 3.9
 $$
Thus we obtain that the system
of the Wilson polynomials
$$
p_n(1-a_1,1-a_2,1-a_3,1-a_4; s^2);
\qquad 4n<a_1+a_2+a_3+a_4-3
\tag 3.10
$$
is orthogonal with respect to the Dougall
weight
 (3.8).

\smallskip

{\bf 3.5. \bf Orthogonality relations
associated with the Askey integral (0.5).}
Now consider the Askey weight, i.e., the weight
that is the integrand in
 (0.5)
$$
\frac{\Gamma(1-2s)\Gamma(1+2s)}
{\prod_{j=1}^4\Gamma(a_j+s)\Gamma(a_j-s)}
$$
In notation of Lemma
 3.2,
 $$ \ell(g)= \int_{-\infty}^\infty g(s)w_3(s)ds,$$
$$h_k(s)=(1-a_1+s)_k(1-a_1-s)_k.$$
Then
 $$\multline
\ell(h_k)=\frac{1}
{\Gamma(a_2+a_3)\Gamma(a_2+a_4)\Gamma(a_3+a_4)}
\times\\ \times
\frac{\Gamma(a_1+a_2 +a_3+a_4-k-3)}
{\Gamma(a_1+a_2-k-1)\Gamma(a_1+a_3-k-1)\Gamma(a_1+a_4-k-1)}
\endmultline
$$
and we obtain the expression coinciding with
 (3.9) up to a constant factor.
Hence, we obtain the same system  (3.10) of orthogonal polynomials.

 \medskip

{\bf \S 4.  Examples of index integrals}

\medskip

 Index integrals are by themselves fairly known subject,
see, for instance \cite{Sla}.
In this section, we present several amusing index integrals
extending the integrals
 (0.1), (0.2), (2.2)--(2.5).

\smallskip

{\bf 4.1.} The index transform of the function
$(1+x)^{b-a} (x+y+1)^{-c-b}$ was evaluated in Lemma 1.1.
Using the inversion formula,
we get
$$\multline
\int_0^\infty
\left| \frac{\Gamma(a+is)\Gamma(b+is)\Gamma(c+is)}
            {\Gamma(2is)}\right|^2
 \F\left[\matrix c+is,c-is\\a+c\endmatrix;-y\right]
\F\left[\matrix b+is,b-is\\a+b\endmatrix;-x\right]
\,ds = \\=
\frac {\pi \Gamma(a+b)\Gamma(a+c)\Gamma(b+c)}
      {(1+x+y)^{c+b}}
.\endmultline
$$

{\bf 4.2.} Writing out the Plancherel formula for the functions
 $(1+x)^{b-a} (x+y+1)^{-c-b}$
and $(1+x)^{-a-d}$,
we obtain (see de Branges \cite{dB}, Theorem 13)
$$\multline
\int_0^\infty
\left| \frac{\Gamma(a+is)\Gamma(b+is)\Gamma(c+is)\Gamma(d+is)}
            {\Gamma(2is)}\right|^2
  \F\left[\matrix c+is,c-is\\a+c\endmatrix;-y\right]ds=
\\
=\frac{\pi\Gamma(a+b)\Gamma(a+c)\Gamma(a+d)\Gamma(b+c)
\Gamma(b+d)\Gamma(c+d)}{\Gamma(a+b+c+d)}
\F \left[\matrix b+c,c+d\\a+b+c+d\endmatrix;-y\right]
.$$
 \endmultline
$$
This formula coincides (up to a permutation of letters)
with the inversion formula for
$\F \left[\matrix a+c,a+d\\a+b+c+d\endmatrix;-x\right]$,
see Lemma 1.2.{\rm c}
(this implies the Second Barnes Lemma
(0.3); evidently, it is not the most simple
 its proof).

\smallskip

{\bf 4.3.} Applying the Plancherel formula to the pair
of functions
$ (1+x)^{b-a}(1+x+y)^{-b-c}$  and
 $ (1+x)^{b-a}(1+x+z)^{-b-d}$,
we obtain
$$\multline
\frac 1\pi
 \int_0^\infty
\left| \frac{\Gamma(a+is)\Gamma(b+is)\Gamma(c+is)\Gamma(d+is)}
            {\Gamma(2is)}\right|^2
\times\\ \times
\F\left[\matrix c-is,c+is\\a+c\endmatrix; -y\right]
 \F\left[\matrix d-is,d+is\\a+d\endmatrix; -z\right] \,ds
=\endmultline
$$
$$
= \Gamma(c+a)\Gamma(c+b)\Gamma(d+a)\Gamma(d+b)
\int_0^\infty \frac {x^{a+b-1} (1+x)^{b-a}dx}
  {(x+y+1)^{b+c}(x+z+1)^{b+d}}
.$$
In the right-hand side we have one of integral representations
of the Appel function
 $F_1$.
Substituting $z=y$, we obtain
$$
\multline
 \frac 1\pi \int_0^\infty
\left| \frac{\Gamma(a+is)\Gamma(b+is)\Gamma(c+is)\Gamma(d+is)}
            {\Gamma(2is)}\right|^2
\times\\ \times
\F\left[\matrix c-is,c+is\\a+c\endmatrix; -y\right]
 \F\left[\matrix d-is,d+is\\a+d\endmatrix; -y\right] \,ds
=\endmultline
$$
$$
 =\frac{\pi\Gamma(a+b)\Gamma(a+c)\Gamma(a+d)\Gamma(b+c)
\Gamma(b+d)\Gamma(c+d)}{\Gamma(a+b+c+d)}
\F\left[\matrix 2b+c+d,c+d\\a+b+c+d\endmatrix;-y
\right]
 .$$

{\bf 4.4.} Applying the Plancherel formula
to the pair of functions
$(1+x)^{b-a}(1+x+y)^{-e-b}$
and $\F\left[\matrix a+c,a+d\\a+b+c+d\endmatrix;-x\right]$,
we get
$$\multline
\frac 1\pi  \int_0^\infty
\left| \frac{\Gamma(a+is)\Gamma(b+is)
\Gamma(c+is)\Gamma(d+is)\Gamma(e+is)}
            {\Gamma(2is)}\right|^2
\F\left[\matrix a+is, a-is\\a+c\endmatrix;-y\right] \,ds
 =\\=                \frac 1{2\pi i}
(1+y)^{a-e}
\Gamma(b+c)\Gamma(b+d)\Gamma(c+d)
\times\\\times
\inti\Gamma\bmatrix a+b+s, a+c+s,a+d+s, e-a-s,-s\\
a+b+c+d+s\endbmatrix (1+y)^s ds
\endmultline
$$
(it is a sum of two functions $\FF$).

\smallskip

{\bf 4.5.}   It is easy to extend this list.
Even our Lemmas 1.1--1.2 were not completely
utilized.

\medskip

\sc
Math.Physics Group

Institute of  Theoretical and Experimantal Physics

Bol\mz shaya Cheremushkinskaya, 25

Moscow 117259

Russia

\tt e-mail neretin\@main.mccme.rssi.ru

\enddocument